\documentclass[twoside,12pt]{article}

\usepackage[all]{xy}
\usepackage{amssymb,amsmath,latexsym}
\usepackage{theorem}
\usepackage{amsfonts}
\usepackage{ulem}
\usepackage{amscd}
\usepackage{amsxtra}
\usepackage{mathrsfs}

\newtheorem{theorem}{Theorem}[section] 
\newtheorem{lemma}[theorem]{Lemma} 
\newtheorem{conjecture}[theorem]{Conjecture}

\numberwithin{equation}{section}

\newcommand{\cqfd}
{\hspace{1cm}
\rule{2mm}{2mm}%
\medbreak%
\par%
}

\def\pr{{\parindent0pt {\bf Proof.\ }}}
\def\cqfd
{\hspace{1cm}
\rule{2mm}{2mm}%
\medbreak%
\par%
}
 
\author{}

\begin{document}
\title{On Generalized $(m, n)$-Jordan Derivations and Centralizers of Semiprime Rings}

\date{}
\maketitle \vspace*{-1.5cm}

\thispagestyle{empty}


\begin{center}
\author{Driss Bennis$^{1,\textbf{a}}$, Basudeb Dhara$^{2}$  and Brahim Fahid$^{1,\textbf{b}}$}
\bigskip

\small{1:  Centre de Recherche de Math\'ematiques et Applications de Rabat (CeReMAR), Faculty of Sciences, Mohammed V University in Rabat, Rabat, Morocco}
\\
  \small{\textbf{a}: d.bennis@fsr.ac.ma; driss$\_$bennis@hotmail.com}\\
  \small{\textbf{b}: fahid.brahim@yahoo.fr}\\

\small{2:   Department of Mathematics, Belda College, Belda,\\ Paschim Medinipur, 721424, W.B., India.}\\
  \small{basu$\_$dhara@yahoo.com} 
\end{center}
\bigskip\bigskip

\noindent{\large\bf Abstract.} In this paper we give an affirmative answer to two  conjectures on generalized  $(m,n)$-Jordan derivations and  generalized  $(m,n)$-Jordan centralizers raised in [S. Ali and A. Fo\v{s}ner, \textit{On Generalized $(m, n)$-Derivations and Generalized $(m, n)$-Jordan Derivations in Rings,} Algebra  Colloq. \textbf{21} (2014), 411--420] and [A. Fo\v{s}ner, \textit{A note  on generalized  $(m,n)$-Jordan centralizers,} Demonstratio Math. \textbf{46} (2013), 254--262].  Precisely, when $R$   is a  semiprime ring,   we prove,  under some suitable torsion restrictions, that   every  nonzero generalized  $(m,n)$-Jordan derivation (resp., a generalized  $(m,n)$-Jordan centralizer) is a derivation  (resp., a two-sided centralizer).\bigskip



\section{Introduction} 
Throughout this paper, $R$ will represent an associative ring with   center $Z(R)$.
 We denote by $char(R)$ the characteristic of a prime ring $R$.  Let $n\geq 2$ be an integer. A ring $R$ is said to be $n$-torsion free if, for all $x \in R$, $nx = 0$ implies $x = 0$. Recall that a ring $R$ is prime if, for any $a,b \in R$, $aRb = \{0\}$ implies $a = 0$ or $b = 0$. A ring $R$ is called
semiprime if, for any $a\in R$, $aRa = \{0\}$ implies $a = 0$.\medskip

An additive mapping $d : R\longrightarrow R$ is called a derivation, if $d(xy) = d(x)y + xd(y)$ holds for all $x, y \in R$, and it is called a Jordan derivation, if $d(x^2) = d(x)x + xd(x)$ holds for all $x \in R$. An additive mapping $T : R \longrightarrow R$ is called a left (resp., right) centralizer if $T(xy) = T(x)y$  (resp., $T(xy) = xT(y)$) is fulfilled for all  $x, y \in R$, and it is called a left (resp., right) Jordan centralizer if $T(x^2) = T(x)x$  (resp., $T(x^2) = xT(x)$) is fulfilled  for all $x \in R$. We call an additive mapping $T : R \longrightarrow R$ a two-sided centralizer (resp., a two-sided Jordan centralizer) if $T$ is both a left as well as a right centralizer (resp., a left and a right Jordan centralizer).\medskip

An additive mapping $ F : R \longrightarrow R$ is called a generalized derivation if $ F(xy) =F(x)y + xd(y)$ holds for all $ x, y \in R$, where $d : R \longrightarrow R $ is a derivation. The concept of generalized derivations was introduced by Bre\v{s}ar in \cite{B} and covers both the concepts of derivations and left centralizers. It is easy to see that generalized derivations are exactly those additive mappings $F$ which can be written in the form $F = d + T$, where $d$ is a derivation and $T$ is a left centralizer.\medskip

The Jordan counterpart of the notion of generalized derivation was introduced by Jing and  Lu in \cite{JL} as follows: An additive mapping $ F : R \longrightarrow R$ is called a generalized Jordan derivation if $F(x^2) = F(x)x + xd(x)$ is fulfilled  for all $x \in R$, where $d : R \longrightarrow R $ is a Jordan derivation.\medskip

The study of relations between various sorts of derivations goes back to Herstein's classical result \cite{H57}  which shows  that any Jordan derivation on a $2$-torsion free prime ring is a derivation (see also \cite{BV88} for a brief proof of Herstein's result). In \cite{C75}, Cusack generalized Herstein's result to $2$-torsion free semiprime rings (see also \cite{B88} for an alternative proof). Motivated by these classical results,  Vukman \cite{V} proved that any generalized Jordan derivation on a $2$-torsion free semiprime ring is a generalized derivation.\medskip

In the last few years several authors have introduced and studied  various sorts of parameterized derivations. In \cite{AF}, Ali and  Fo\v{s}ner defined the notion of $(m, n)$-derivations as follows:   Let $m,n\geq 0 $ be two fixed integers with $ m+n\neq 0$. An additive mapping $ d: R\longrightarrow R$ is called an $(m,n)$-derivation if
  \begin{equation}
    (m+n)d(xy)=2md(x)y+2nxd(y)
  \end{equation}
holds for all $x, y \in R$.\medskip

Obviously, a $(1, 1)$-derivation on a $2$-torsion free ring is a derivation.\medskip

In the same paper \cite{AF}, a generalized $(m, n)$-derivation was defined as follows:   Let $m,n\geq 0$  be two fixed integers with $ m+n\neq 0$. An additive mapping $ D: R\longrightarrow R$ is called a  generalized $(m,n)$-derivation if  there exists an $(m, n)$-derivation $d : R \longrightarrow R$ such that
  \begin{equation}
    (m+n)D(xy)=2mD(x)y+2nxd(y)
  \end{equation}
holds for all $x, y \in R$.\medskip

Obviously, every generalized $(1, 1)$-derivation on a $2$-torsion free ring is   a generalized derivation.\medskip

In \cite{V1}, Vukman defined an $(m, n)$-Jordan derivation as follows:   Let $m,n\geq 0 $  be two fixed integers with $ m+n\neq 0$. An additive mapping $ d: R\longrightarrow R$ is called an $(m,n)$-Jordan derivation if
  \begin{equation}
    (m+n)d(x^2)=2md(x)x+2nxd(x)
  \end{equation}
holds   for all $x, y \in R$.\medskip

Clearly, every  $(1,1)$-Jordan derivation on a $2$-torsion free ring is  a Jordan derivation.\medskip

Recently, in \cite{KV},  Kosi-Ulbl and Vukman proved the following result. 

\begin{theorem}[\cite{KV},  Theorem 1.5]\label{th-lem}
  Let $ m, n \geq 1$ be distinct integers, $R$ a  $mn(m + n)|m -n| $-torsion free semiprime ring  an  $ d  : R\longrightarrow R$ an $(m, n)$-Jordan derivation. Then $d$ is a derivation which maps $R$ into $Z(R)$.
\end{theorem}

The $(m, n)$-generalized counterpart of the notion of an $(m, n)$-Jordan derivation is introduced by  Ali and  Fo\v{s}ner  in \cite{AF} as follows:   Let $m,n\geq 0 $  be two fixed integers with $ m+n\neq 0$. An additive mapping $ F: R\longrightarrow R$ is called a  generalized $(m,n)$-Jordan derivation if  there exists an $(m, n)$-Jordan derivation $d : R \longrightarrow R$ such that
  \begin{equation}
    (m+n)F(x^2)=2mF(x)x+2nxd(x)
  \end{equation}
holds  for all $x, y \in R$.\medskip

Based on some observations and inspired by the classical results, Ali and  Fo\v{s}ner in  \cite{AF}  made the following conjecture.

\begin{conjecture}[\cite{AF}, Conjecture 1]
Let $ m, n \geq 1$ be two fixed integers,  let $R$ be a semiprime ring with suitable torsion restrictions, and  let $F : R\longrightarrow R$ be a nonzero generalized $(m, n)$-Jordan derivation. Then $F$ is a derivation which maps $R$ into $Z(R)$.
\end{conjecture}

The first aim of this paper is to give an affirmative answer to this conjecture. Namely, our first main result is the following theorem.

 \begin{theorem}\label{th-princ1}
Let $m, n \geq 1$ be distinct integers, let $R$ be a $k$-torsion free semiprime ring, where $k=6mn(m+n)| m-n|$, and let $F: R\longrightarrow R$ be a nonzero generalized  $(m,n)$-Jordan derivation.  Then $F$ is a derivation which maps $R$ into $Z(R)$.
\end{theorem}

On the other hand and in parallel, there are similar works which study relations between various sorts of Jordan centralizers and centralizers. Namely,  in \cite{Z91}, Zalar   proved   that any left (resp., right) Jordan centralizer on a $2$-torsion free semiprime ring is a left (resp.,  right) centralizer. In \cite{V99}, Vukman proved that, for    a $2$-torsion free semiprime ring  $R$, every additive mapping $ T: R\longrightarrow R$  satisfying the relation  ``$2T(x^2) = T(x)x + xT(x)$ for all  $x \in R$" is a two-sided centralizer. Motivated by these results and inspired by his work \cite{V99}, Vukman  in \cite{V10} introduced the notion of an  $(m, n)$-Jordan centralizer as follows:    Let $m, n\geq 0 $ be two  fixed integers with $ m+n\neq 0$. An additive mapping $ T: R\longrightarrow R$ is called an $(m,n)$-Jordan centralizer if
  \begin{equation}
    (m+n)T(x^2)=mT(x)x+nxT(x)
  \end{equation}
holds for all $x, y \in R$.\medskip

Obviously, a $(1, 0)$-Jordan centralizer (resp., $(0, 1)$-Jordan centralizer) is a left (resp., a right) Jordan centralizer. When $n=m=1$, we recover the maps studied in \cite{V99}.\medskip

Based on some observations and results, Vukman conjectured that, on semiprime rings with suitable torsion restrictions, every $(m, n)$-Jordan centralizer is a two-sided centralizer (see \cite{V10}).  Recently, this conjecture was solved affirmatively by  Kosi-Ulbl and Vukman   in \cite{KV2}. Namely, they   proved the following result.

\begin{theorem}[\cite{KV2},  Theorem 1.5]\label{th-lem2}
  Let $m, n \geq 1$ be distinct integers, let $R$ be an $mn(m + n) $-torsion free semiprime ring, and let $T: R\longrightarrow R$ be an  $(m,n)$-Jordan  centralizer.  Then $T$ is a two-sided centralizer.
\end{theorem}

Inspired by the work of Vukman    \cite{V99,V10},  Fo\v{s}ner     \cite{F13}   introduced  more generalized version  of $(m, n)$-Jordan centralizers as follows: Let $m,n\geq 0 $  be two fixed integers with $ m+n\neq 0$. An additive mapping $ T: R\longrightarrow R$ is called a generalized $(m,n)$-Jordan centralizer if there exists an $(m,n)$-Jordan centralizer $T_0: R\longrightarrow R$ such that
  \begin{equation}
    (m+n)T(x^2)=mT(x)x+nxT_0(x)
  \end{equation}
holds    for all $x \in R$.\medskip

Thus, a generalized $(1, 0)$-Jordan centralizer is a left Jordan centralizer.\medskip

In \cite{F13}, Fo\v{s}ner showed that,    on  a prime ring  with a specific  torsion condition, every generalized $(m, n)$-Jordan centralizer is a two-sided centralizer. This led  Fo\v{s}ner to  make  the following conjecture.

\begin{conjecture}[\cite{F13}, Conjecture 1]
Let $m, n \geq 1$ be two fixed integers, let $R$ be a semiprime ring with suitable torsion restrictions, and let $T : R\longrightarrow R $ be a generalized $(m, n)$-Jordan centralizer. Then $T$ is a two-sided centralizer.
\end{conjecture}

The second aim of this paper is to give an affirmative answer to  Fo\v{s}ner's conjecture. Namely, our second main result is the following theorem.

 \begin{theorem}\label{th-princ2}
  Let $m, n \geq 1$ be two fixed integers, let $R$ be an $6mn(m+n)(2n+m)$-torsion free semiprime ring, and let $T : R\longrightarrow R$ be a nonzero generalized $(m,n)$-Jordan centralizer.  Then $T$ is a two-sided centralizer.
\end{theorem}

\section{Proof of the main theorems}

In the proof of our main results, Theorems \ref{th-princ1} and \ref{th-princ2}, we shall use the following results.

\begin{lemma}[\cite{AF}, Lemma 1]\label{th-lem1}
 Let $m, n \geq 0$ be distinct integers with $m + n\neq 0$,  let $R$ be a $2$-torsion free ring, and let $F : R\longrightarrow R$ be a nonzero generalized  $(m,n)$-Jordan derivation  with  an associated $(m, n)$-Jordan derivation $d$.  Then,  $(m+n)^{2}F(xyx) = m(n-m)F(x)xy+m(3m+n)F(x)yx+m(m-n)F(y)x^2 +4mnxd(y)x  + n(n-m)x^{2}d(y)+n(m+3n)xyd(x)+n(m-n)yxd(x)$ for all $x,y \in R$.
\end{lemma}

\begin{lemma}[\cite{BE}, Theorem 3.3]\label{Lemma A}
  Let $ n \geq 2$ be a fixed integer and let $R$ be a prime ring with $char (R) = 0$ or $char (R) \geq n$. If $T : R\longrightarrow R$ is an additive mapping satisfying the relation $T(x^n) = T(x)x^{n-1}$ for all $x\in R$, then $T(xy) = T(x)y$ for all  $ x, y \in R$.
\end{lemma}

\begin{lemma}[\cite{F13}, Lemma 1]\label{lem-th2}
 Let $m, n \geq 0$ be distinct integers with $m + n\neq 0$, let $R$ be a ring, and let $T : R\longrightarrow R$ be a nonzero generalized  $(m,n)$-Jordan centralizer with an associated $(m,n)$-Jordan centralizer  $T_0$.  Then, $   2(m+n)^{2}T(xyx) = mnT(x)xy+m(2m+n)T(x)yx-mnT(y)x^2+2mnxT_0(y)x - mnx^{2}T_0(y)+n(m+2n)xyT_0(x)+mnyxT_0(x)$
                    for all $x,y \in R$.
\end{lemma}

\begin{lemma}[\cite{V2},  Lemma 3]\label{lem-v}
  Let $R$ be a semiprime ring and let $T : R \longrightarrow R$ be an additive mapping. If either
$T(x)x = 0$ or $x T(x) = 0$ holds for all $x \in R$, then $T = 0$.
\end{lemma}

We shall use the relation between semiprime rings and prime ideals.  Namely, it is well-known that a ring $R$ is semiprime if and only if   the intersection of all prime ideals of $R$ is zero   if and only if  $ R$ has no nonzero nilpotent (left, right) ideals (see for instance Lam's book \cite{Lam} or the recent   book of Bre\v{s}ar \cite{Bboo}). Due to the classical Levitzki's paper \cite{L51}, several authors prefer to refer to a such result by Levitzki's lemma.\\

Let $I$ be an ideal of $R$. For an element $x\in R$, we use $\overline{x}$ to denote the equivalence class of $x$  modulo $I$.

\begin{lemma}\label{Lemma B}
    Let $R$ be both a $2$-torsion free and a $3$-torsion free semiprime ring and let $T : R \longrightarrow R$ be an additive map such that $T(x)x^3 = 0$ and $T(x^4) = 0$ for all $x \in R$. Then $T(xy) = T(x)y$ for all $x, y \in R$.
\end{lemma}
\pr
Let $x,y \in R$. We prove that $T(xy)=T(x)y$. We may assume that $x$ and $y$ are not $0$. Let $P$ be a prime ideal of $R$ and set $\overline{R} = R/P$. Consider an element $p \in P$. By hypothesis, $0 = T(x+p)(x+p)^3 = (T(x)+T(p))(x^3+xpx+px^2+p^2x+x^2p+xp^2+pxp+p^3)$. Thus, $ 0 = T(x)(xpx+px^2 +p^2x+x^2p+xp^2 +pxp+p^3)+
T(p)(x^3 +xpx+px^2 +p^2x+x^2p+xp^2 +pxp)$. Hence, $T(p)x^3 \in P$, equivalently $\overline{T(p)}\overline{x}^3 = 0$.  By  Levitzki's lemma, $\overline{T(p)}\overline{x} = 0$, and then $\overline{T(p)} = 0$ (since $\overline{R}$ is a prime ring). Thus, $T(P)\subseteq P$, which implies that $ T(x+P) = T(x)+P$. Then, the induced map $\overline{T}:R/P \rightarrow R/P$ such that $\overline{T}(\overline{x})=\overline{T(x)}$ for every $x\in R$, is well defined. Now, since $\overline{T}(\overline{x})\overline{x}^3 = 0$ and $\overline{T}(\overline{x}^4) = 0$, $\overline{T}(\overline{x}^4) = \overline{T}(\overline{x})\overline{x}^3$. This shows, using Lemma \ref{Lemma A}, that $\overline{T}(\overline{xy}) = \overline{T}(\overline{x})\overline{y}$. Therefore,  $T(xy)-T(x)y \in P$. Finally, by the semiprimeness of $R$, we get the desired result.\cqfd

Now we are ready to prove the first main result.\medskip

\noindent\textbf{Proof of Theorem \ref{th-princ1}.}  Let $d$ be the associated $(m,n)$-Jordan derivation of $F$. Since $R$ is a semiprime ring, $d$ is a derivation which maps $R$ into $Z(R)$ (by Theorem \ref{th-lem}).  Let us denote $F -d$ by $D$. Then,  we have $(m+n)D(x^2)=(m+n)F(x^2)-(m+n)d(x^2)=2mF(x)x+2nxd(x)-2md(x)x-2nxd(x)=2mD(x)x$  for all $x\in R$. Thus
\begin{equation}\label{eq-3}
 (m+n)D(x^2)=2mD(x)x,\ x\in R.
\end{equation}
Replacing $x$ with $x^2$ in (\ref{eq-3}), we get
\begin{equation}\label{eq-4}
 (m+n)D(x^4)=2mD(x^2)x^2,\ x\in R.
\end{equation}
 Multiplying by $m + n$ and then using (\ref{eq-3}), we get
\begin{equation}\label{eq-5}
  (m + n)^2D(x^4) = 4m^2D(x)x^3,\ x\in R.
\end{equation}
 On the other hand,   putting $x^2$ for $y$ in the relation
of Lemma \ref{th-lem1} and using the fact that  $D$ is a generalized $(m, n)$-Jordan derivation associated with the zero map as an $(m, n)$-Jordan derivation,  we get
\begin{equation}\label{eq-6}
    (m+n)^{2}D(x^4) = m(n-m)D(x)x^3+m(3m+n)D(x)x^3+m(m-n)D(x^2)x^2,\ x\in R.
\end{equation}
Multiplying both sides in (\ref{eq-6}) by $2$  we get
\begin{equation}\label{eq-7}
    2(m+n)^{2}D(x^4) = 2m(n-m)D(x)x^3+2m(3m+n)D(x)x^3+2m(m-n)D(x^2)x^2,\ x\in R.
\end{equation}
Combining (\ref{eq-4}) and (\ref{eq-7}), we get
\begin{equation}\label{eq-8}
    2(m+n)^{2}D(x^4) = 2m(n-m)D(x)x^3+2m(3m+n)D(x)x^3+(m+n)(m-n)D(x^4),\ x\in R,
\end{equation}
which gives
\begin{equation}\label{eq-9f}
    (m+n)(m+3n)D(x^4) = 4m(m+n)D(x)x^3,\ x\in R.
\end{equation}
Multiplying both sides in (\ref{eq-9f}) by $m+n$, we get
\begin{equation}\label{eq-9}
    (m+n)^2(m+3n)D(x^4) = 4m(m+n)^2D(x)x^3,\ x\in R.
\end{equation}
Multiplying by $m+3n$ in (\ref{eq-5}), we get
\begin{equation}\label{eq-10}
  (m + n)^2(m+3n)D(x^4) = 4m^2(m+3n)D(x)x^3,\ x\in R.
\end{equation}
By comparing (\ref{eq-9}) and (\ref{eq-10}), we get
\begin{equation}\label{eq-11}
  4mn(m-n)D(x)x^3=0,\ x\in R.
\end{equation}
 Since $R$ is a $2mn|n-m|$-torsion free ring,  $D(x)x^3=0$ for all $x\in R$.  Applying $D(x)x^3 = 0$  in equation (\ref{eq-5}), we get $(m + n)^2D(x^4) = 0$ for all $x\in R$. By using the torsion free restriction, we have $D(x^4)=0$ for all $x\in R$. Hence, $D(xy) = D(x)y$ for all $x, y \in R$ (by Lemma \ref{Lemma B}). Applying this in (\ref{eq-3}), yields $(m+n)D(x)x = 2mD(x)x$ for all $x\in R$, equivalently $(m-n)D(x)x = 0$. Since $R$ is an $|m-n|$-torsion free ring, $D(x)x = 0$ for all $x\in R$. Therefore, by Lemma \ref{lem-v},  $D =0$.  This completes the proof.
\cqfd

The second main result is proved similarly. Nevertheless,  we include a proof for completeness.\\

\noindent\textbf{Proof of Theorem \ref{th-princ2}.} Let $T_0$ be the associated $(m,n)$-Jordan centralizer of $T$. Since $R$ is a semiprime ring, $T_0$ is a two-sided centralizer (by Theorem \ref{th-lem2}).  Let us denote $T -T_0$ by $D$. Then,  we have $(m+n)D(x^2)=(m+n)T(x^2)-(m+n)T_0(x^2)=mT(x)x+nxT_0(x)-m T_0(x)x-nxT_0(x)=mD(x)x$ for all $ x\in R.$  Thus
\begin{equation}\label{eq-13}
 (m+n)D(x^2)=mD(x)x,\ x\in R.
\end{equation}
Replacing $x$ with $x^2$ in (\ref{eq-13}), we get
\begin{equation}\label{eq-14}
 (m+n)D(x^4)=mD(x^2)x^2,\ x\in R.
\end{equation}
 Multiplying by $m + n$ and then using (\ref{eq-13}), we get
\begin{equation}\label{eq-15}
  (m + n)^2D(x^4) = m^2D(x)x^3,\ x\in R.
\end{equation}
 On the other hand, if we put $y = x^2$ in the relation
of Lemma \ref{lem-th2}, we get
\begin{equation}\label{eq-16}
    2(m+n)^{2}D(x^4) = mnD(x)x^3+m(2m+n)D(x)x^3-mnD(x^2)x^2,\ x\in R.
\end{equation}
Multiplying  both sides in (\ref{eq-15})  by $2$  we get
\begin{equation}\label{eq-17}
   2(m + n)^2D(x^4) = 2m^2D(x)x^3,\ x\in R.
\end{equation}
Combining (\ref{eq-14}) and (\ref{eq-16}), we get
\begin{equation}\label{eq-18}
    2(m+n)^{2}D(x^4) = mnD(x)x^3+m(2m+n)D(x)x^3-n(m+n)D(x^4),\ x\in R,
\end{equation}
which implies
\begin{equation}\label{eq-19x}
    (m+n)(2m+3n)D(x^4) = 2m(m+n)D(x)x^3,\ x\in R.
\end{equation}
Multiplying both sides of above relation by $m+n$, we have
\begin{equation}\label{eq-19}
    (m+n)^2(2m+3n)D(x^4) = 2m(m+n)^2D(x)x^3,\ x\in R.
\end{equation}
Multiplying by $(2m+3n)$ in (\ref{eq-15}), we get
\begin{equation}\label{eq-20}
  (m + n)^2(2m+3n)D(x^4) = m^2(2m+3n)D(x)x^3,\ x\in R.
\end{equation}  By comparing (\ref{eq-19}) and (\ref{eq-20}), we get
\begin{equation}\label{eq-21}
  mn(2n+m)D(x)x^3=0,\ x\in R.
\end{equation}
Since $R$ is a $mn(2n+m)$-torsion free ring,  $D(x)x^3=0$ for all $x\in R$.  Applying $D(x)x^3 = 0$  in equation (\ref{eq-15}) and then using $(m+n)$-torsion freeness of $R$, we get $D(x^4) = 0$ for all $x\in R$. Moreover, since $R$ is a $2$ and a $3$-torsion free ring,  by Lemma \ref{Lemma B}, we get $D(xy) = D(x)y$ for all $x, y \in R$. Applying this in (\ref{eq-13}), yields $(m+n)D(x)x = mD(x)x$ for all $x\in R$. So $nD(x)x = 0$, which implies that $D(x)x = 0$ for all $x\in R$. Therefore,  by Lemma \ref{lem-v},  $D =0$.  This completes the proof.\cqfd\medskip

\noindent\textbf{Acknowledgements.}  The authors would like to thank Professor  Abdellah Mamouni for useful discussions.\\


\begin{thebibliography}{XX}


\bibitem{AF}{S. Ali and A. Fo\v{s}ner,} \textit{On Generalized $(m, n)$-Derivations and Generalized $(m, n)$-Jordan Derivations in Rings,} Algebra  Colloq. \textbf{21} (2014), 411--420.



\bibitem{B88}{M. Bre\v{s}ar,} \textit{Jordan derivations on semiprime rings,} Proc. Amer. Math. Soc.  \textbf{104} (1988), 1003--1006.

\bibitem{B}{M. Bre\v{s}ar,} \textit{On the distance of the composition of two derivations to the generalized derivations,} Glasg. Math. J. \textbf{33} (1991), 89--93.

\bibitem{Bboo}{M. Bre\v{s}ar,}  \textit{Introduction to noncommutative algebra,}  Universitext, Springer, 2014.

\bibitem{BV88}{M. Bre\v{s}ar and J. Vukman,} \textit{Jordan derivations on prime rings,} Bull. Aust. Math. Soc.  \textbf{37} (1988), 321--322.


\bibitem{BE}{D. Benkovi\v{c} and D. Eremita}, \textit{Characterizing left centralizers by their action on a polynomial,}  Publ. Math. Debercen \textbf{64} (2004), 343--351.

\bibitem{C75}{J. Cusak,}  \textit{Jordan derivations on rings,}  Proc. Amer. Math. Soc.  \textbf{53} (1975), 321--324.

\bibitem{F13}{A. Fo\v{s}ner,} \textit{A note  on generalized  $(m,n)$-Jordan centralizers,} Demonstratio Math. \textbf{46}
(2013), 254--262.



\bibitem{H57}{I. N. Herstein},  \textit{Jordan derivations of prime rings,}  Proc. Amer. Math. Soc.  \textbf{8} (1957), 1104--1119.


\bibitem{JL}{W. Jing and S. Lu,} \textit{Generalized Jordan derivations on prime rings and standard operator algebras,} Taiwanese J. Math. \textbf{7} (2003), 605--613.


\bibitem{KV}{I. Kosi-Ulbl  and J. Vukman,} \textit{A note on $(m,n)$-Jordan derivation of rings and banach algebras,}  Bull. Aust. Math. Soc. \textbf{93} (2016), 231--237.

\bibitem{KV2}{I. Kosi-Ulbl and J.Vukman,} \textit{On $(m,n)$-Jordan centralizers of semiprime rings,}  Publ. Math. Debrecen    \textbf{7490} (2016), 1--9.


\bibitem{Lam}{T. Y. Lam,} \textit{A first course in noncommutative rings,}  Graduate Texts in Mathematics, 123, Springer-Verlag, New York, 1991.

\bibitem{L51}{J. Levitzki,} \textit{Prime ideals and the lower radical,}    Amer. J. Math.   \textbf{73} (1951), 25--29.

\bibitem{V99}{J. Vukman,} \textit{An identity related to centralizers in semiprime rings,} Comment. Math. Univ. Carolin. \textbf{40} (1999), 447--456.


 \bibitem{V2}{J. Vukman,} \textit{Identities with derivations and automorphisms on semiprime rings,} Int. J. Math. Math. Sci.   \textbf{7} (2005), 1031--1038.

\bibitem{V}{J. Vukman,} \textit{A note on generalized derivations of semiprime rings,} Taiwanese J. Math. \textbf{11} (2007), 367--370.


\bibitem{V1}{J. Vukman,} \textit{On $(m, n)$-Jordan derivations and commutativity of prime rings,}  Demonstratio Math. \textbf{41} (2008), 773--778.


\bibitem{V10}{J. Vukman,} \textit{On $(m, n)$-Jordan centralizers in rings and algebras,}  Glasg. Math. J. \textbf{45}  (2010), 43--53.

\bibitem{Z91}{B. Zalar,} \textit{On centralizers of semiprime rings,}  Comment. Math. Univ. Carolin.  \textbf{32} (1991),
609--614.




 \end{thebibliography}
\end{document}